\documentclass[12pt,a4paper]{amsart}
\usepackage{amsthm, amssymb}
\usepackage{amsfonts}
\usepackage{amsthm}
\usepackage{amsmath}
\usepackage{amscd}
\usepackage[latin2]{inputenc}
\usepackage{t1enc}
\usepackage[mathscr]{eucal}

\numberwithin{equation}{section}
\usepackage[margin=2.9cm]{geometry}
\usepackage{epstopdf} 
\usepackage{hyperref}

 \newtheorem{definition}{Definition}[section]

\theoremstyle{definition}

\theoremstyle{plain}
\newtheorem{Th}{Theorem}[section]
\newtheorem{lemma}{Lemma}[section]

\newtheorem{Prop}[Th]{Proposition}
 \theoremstyle{definition}
\newtheorem{Def}[Th]{Definition}

\newtheorem{?}[Th]{Problem}

\begin{document}
\begin{center}
\title*{Classification of irreducible integrable modules for extended affine Lie
algebras with center acting trivially }
\end{center}
\begin{center}
\author*{Santanu Tantubay$^1$, Punita Batra$^2$}
\end{center} 
{$^{1,2}$ Harish-Chandra Research Institute, HBNI, Prayagraj (Allahabad) 211019, India\\
$^1${santanutantubay@hri.res.in}\\[2pt]
$^2${batra@hri.res.in}}
\maketitle
\begin{abstract}
 \begin{center}
   We classify the irreducible integrable modules for the twisted toroidal
extended affine Lie algebras with finite diemnsional weight spaces where the finite
dimensional center acts trivially. We have proved that the entire central extension part acts trivially on the modules and our modules turn out to be highest weight modules.
\end{center}
 
{\bf Msc: 17B67, 17B66.}
\end{abstract}
\newpage
\section{Introduction}     
Extended affine Lie algebras (EALAs) form a category of important Lie algebras consisting of finite dimensional semisimple Lie algebras, affine Lie algebras and some other classes of Lie algebras. Twisted toroidal extended affine Lie algebras are examples of EALAs. The structure theory of EALAs have been developed by several mathematicians like Allison, Azam, Berman, Gao, Neher, Pianzola and Yoshii (see \cite{[1]}, \cite{[16]}, \cite{[17]} and references therein). With an extended affine Lie algebra, we associate  a graded ideal, the core and its central quotient, the centerless core, both are Lie tori. So it is important to study the representation theory of EALAs.
\\\\
Let $g$ be a finite dimensional simple Lie algebra over $\mathbb{C}$ and $A_n$ be a  Laurent polynomial ring in $n\geq 2$ commuting variables $t_1, \dots ,t_n$. Let $L(g)=g \otimes A_n$ be a loop algebra. Let $L(g)\oplus \Omega_{A_n}/dA_n$ be the universal central extension of $L(g)$. This is called toroidal Lie algebras (TLA). Let $Der(A_n)$ be the Lie algebra of derivations of $A_n$. Then $L(g)\oplus \Omega_{A_n}/dA_n\oplus Der(A_n)$ is called the full toroidal Lie algebra (FTLA). TLA and FTLA are not EALA, as they do not admit any non-degenerate symmetric invariant bilinear form. So instead of $Der(A_n)$, one takes $S_n$ the subalgebra of $Der(A_n)$, consisting of divergence zero vector fields. The Lie algebra  $L(g)\oplus \Omega_{A_n}/dA_n\oplus S_n$ is called a toroidal extended affine Lie algebra which admits a non-degenerate symmetric bilinear form, and hence is an example of EALAs.\\\\
 Here we consider more general EALAs. We take $\sigma_1, \dots ,\sigma_n$ to be finite order commuting automorphisms
  of $g$ and consider the multi-loop algebra $\oplus_{k \in \mathbb{Z}^n}g(\bar{k})\otimes t^k$. The representations of multiloop algebras have been studied by Lau in \cite{[18]}. We assume that this multi-loop algebra is a Lie torus. Now we consider the universal central extension of multi-loop algebra and add $S_n$, the Lie algebra consisting of divergence zero vector fields on $A_n$. This Lie algebra is called twisted toroidal extended affine Lie algebra, we denote it by $\tau$. In \cite{[2]}, S.E. Rao, S. Sharma, P. Batra classified  irreducible integrable modules for $\tau$ with finite dimensional weight spaces, where the zero degree central operators act nontrivially. In \cite{[15]}, the classification of irreducible integrable modules with finite dimensional weight spaces have been done for $n=2$, both zero and non zero case. In this paper we will classify irreducible integrable modules with finite dimensional weight spaces for $\tau$, where the zero degree central operators act trivially for any $n\geq 3$. We make use of some important results from \cite{[11]} and \cite{[12]} in order to classify our modules.\\\\
  The paper has been organized as follows. In Section \ref{sec 2},  we define the twisted toroidal extended affine Lie algebra $\tau$ as above.This Lie algebra has a natural triangular decomposition given by   $\tau=\tau^-\oplus \tau^0\oplus \tau^+$ (see Section \ref{sec 3}). In Section \ref{sec 3}, we prove the existence of highest weight vector with respect to this triangular decomposition. We show that the highest weight space is irreducible module for $\tau^0$ as well as it is $\Gamma$-graded. In Section \ref{sec 4}, we prove that if zero degree center $\mathbb{C}$-span$\{K_1, \dots K_n\}$ acts trivially on $V$, then the whole central extension part  $Z=\Omega_{A_n(m)}/d_{A_n(m)}$ acts trivially on $V$. In Section \ref{sec 5}, to classify the modules we prove an important result Proposition \ref{prop 5.1}(2). We have used a new method to prove this. In Section \ref{sec 6}, we have proved our main theorem, Theorem \ref{Thm 6.6}  using a result of \cite{[11]}.

\section{Notation and Preliminaries}\label{sec 2}
Let $A_{n}=C[t_1 ^{\pm 1},t_2^{\pm 1}\dots , t_n^{\pm 1}]$ be a Laurent polynomial
ring in n variables. Let  $L({g})={g}  \otimes A_n$ be the corresponding loop algebra, where
$g$ is a finite dimensional simple Lie algebra over $\mathbb{C}$ with a Cartan
subalgebra $h$. Let $\Omega_{A_n}$ be a vector space spanned by the symbols
$t^kK_i,1\leq i\leq n,k \in \mathbb{Z}^n$. Let $dA_n$ be the subspace spanned by
$\sum_{i=1}^{n} k_it^kK_i.$ It is well known that $\tilde{L}(g)=L(g) \oplus
\Omega_{A_n}/dA_n $ is the universal central extension of $L(g)$ with the following
brackets:
 \begin{center}
 
  $[x(p),y(q)]=[x,y](p+q)+(x|y)\sum_{i=1}^{n}p_it^{p+q}K_i$,
  \end{center}   where $x(p)=x\otimes t^p$. The basis for the Lie algebra $Der(A_n)$
is \{$d_i,t^rd_i|1\leq i\leq n,0 \neq r \in \mathbb{Z}^n$\}. Now $Der(A_n)$ acts
on $\Omega_{A_n}/dA_n$ by
 \begin{center}
 $t^pd_a(t^qK_b)=q_at^{p+q}K_b+ \delta_{ab} \sum_{c=1}^{n}p_ct^{p+q}K_c$.
 \end{center}
There are two non-trivial 2-cocyle of $Der(A_n)$ with values in $\Omega_{A_n}/dA_n$:
\begin{center}
$\phi_1(t^pd_a,t^qd_b)=-q_ap_b \sum _{i=1}^{n}p_it^{p+q}K_i$
\end{center}
\begin{center}
$\phi_2(t^pd_a,t^qd_b)=-p_aq_b\sum _{i=1}^{n}p_it^{p+q}K_i$.
\end{center}
Let $\phi$ be any linear combination of $\phi_1$ and $
\phi_2$. Then $\tau = L(g)\oplus \Omega_{A_n}/dA_n\oplus Der(A_n)$ is a Lie algebra
with the following brackets and is called a full toroidal Lie algebra:
\begin{center}
$[t^pd_a,X(q)]=q_aX(p+q),$
\end{center}
\begin{center}
$[t^pd_a,t^qK_b]=q_at^{p+q}K_b+\delta_{ab}\sum_{i=1}^{n}p_it^(p+q)K_i$,
\end{center}
\begin{center}
$[t^pd_a,t^qd_b]=q_at^{p+q}d_b-p_bt^{p+q}d_a+\phi(t^pd_a,t^qd_b).$

\end{center}

Now consider the subalgebra of divergence zero vector fields $ \emph{S}_{n}$ of
$Der(A_n)$. One can define $\emph{S}_n=\{D(u,r)|(u|r)=0,u\in \mathbb{C}^n,r\in
\mathbb{Z}^n\}$. Now consider the subalgebra
$\tau_{div}=L(g)\oplus\Omega_{A_n}/dA_n\oplus\emph{S}_n$ of $\tau$.
It is well known that unlike $\tau$, $\tau_{div}$ possesses a non-degenerate
symmetric, invariant bilinear form and is called as toroidal extended affine Lie
algebra. The form is defined as follows:
\begin{center}
$(X(r)|Y(s))=\delta_{r,-s}(X|Y)$,
\end{center}
\begin{center}

 $\forall X,Y\in \overset{\circ}{g},r,s\in \mathbb{Z}^n;$
\end{center}
\begin{center}
$(\sum_{i=1}^n a_it^rd_i|t^sK_d)=\delta_{r,-s}a_d$.
\end{center}
All other brackets of bilinear form are zero. \\Now Let $g_1$ be any arbitrary
finite dimensional simple Lie algebra over $\mathbb{C}$ with a cartan subalgebra
$h_1$. Let $\Delta(g_1,h_1)=supp_{h_1}(g_1)$. Then
$\Delta_1^{\times}=\Delta^{\times}(g_1,h_1)=\Delta(g_1,h_1) - \{0\}$ is an
irreducible reduced finite root system with atmost two root lengths. Define
\begin{center}
$\Delta_{1,en}^{\times}=\Delta_1^{\times}\cup 2\Delta_{1,sh}^{\times}$ \hspace{2mm}
if $\Delta_1^{\times}=B_l$ types
\end{center}
\begin{center}
$\Delta_1^{\times}$ \hspace{2mm} otherwise.
\end{center}
\begin{definition} 
A finite dimensional $g_1$-module $V$ is said to satisfy condition $(M)$ if $V$ is
irreducible with dimension greater than 1 and weights of $V$ relative to $h_1$ are
contained in $\Delta_{1,en}^{\times}$.
\end{definition} 
Now recall that $g$ is a finite dimensional simple Lie algebra with a Cartan
subalgebra $h$ and let $\sigma_1,\sigma_2,\dots, \sigma_n$ be the commuting
automorphisms of $g$ of order $m_1,m_2,\dots, m_n$ respectively. Let
$m=(m_1,m_2,\dots ,m_n)\in \mathbb{Z}^n$. Define $\Gamma=m_1\mathbb{Z}\oplus \dots
\oplus m_n\mathbb{Z}$ and $\Lambda:=\mathbb{Z}^n/\Gamma$. Then we have $
g=\oplus_{\bar{k}}g(\bar{k}) $, where $\bar{k}\in \Lambda$, and $g(\bar{k})=\{X \in g|\sigma_i(X)=\zeta_{i}^{k_i}X,1\leq i\leq n \}$, where $\zeta_i$
are $m_i$-th primitive root of unity for $i=0,\dots,n$.
\begin{definition}
A multiloop algebra $\oplus_{k\in \mathbb{Z}^n}g (\bar{k})\otimes t^k$ is called a
Lie torus $LT$ if 
\begin{itemize}
\item[(1)] $g(\bar{0})$ is a finite dimensional simple Lie algebra.
\item[(2)] For $\bar{k}\neq 0$ and $g(\bar{k})\neq 0$, $g(\bar{k})\cong
U(\bar{k})\oplus W(\bar{k})$, where $U(\bar{k})$ is trivial as $g(\bar{0})$-module
and either $W(\bar{k})$ is zero or satisfy condition (M).
\item[(3)] The order of the group generated by $\sigma_i,1\leq i \leq n$ is equal to
the product of orders of each $\sigma_i$, for $1\leq i\leq n$.
\end{itemize}
   
\end{definition}
Let $h(\bar{0})$ denote a Cartan subalgebra of $g(\bar{0})$. Then by \cite{[13]} Lemma
3.1.3, $h(\bar{0})$ is ad-diagonalizable on $g$ and $\bigtriangleup
^\times=\bigtriangleup^\times (g,h(\bar{0}))$ is an irreducible finite root systemin
$h(\bar{0})$ (Proposition 3.3.5,\cite{[13]}). Let $\bigtriangleup_0
:=\bigtriangleup(g(\bar{0}),h(\bar{0}))$. One of the main properties of Lie tori is
that $\bigtriangleup :=\bigtriangleup_{0,en} $ (Proposition 3.2.5,\cite{[1]}). Let
$A_n(m)=\mathbb{C}[t_1^{\pm 1},\dots,t_n^{\pm 1} ]$ and
$\emph{S}_n(m)=\{D(u,r)|(u|r)=0, u \in \mathbb{C}^n,r \in \Gamma $. Now take $\tau
=LT\oplus \Omega_{A_n(m)}/d_{A_n(m)}\oplus \emph{S}_n(m)$. One can easily check that if
$x\in g(\bar{k})$, $y\in g(\bar{l})$ with $(x|y)\neq 0$, then $\bar{k}+\bar{l} \in
\Gamma$.Therefore we can check that $\tau$ is a well defined Lie algebra called as
twisted toroidal extended affine Lie algebra. The aim of this paper is to classify
the irreducible integrable modules of twisted toroidal EALAs where $K_i$'s are
acting trivially $\forall i \in \{1,\dots,n\}$.

\section{Existence of Highest weight space}\label{sec 3}
In this section we will give a root space decomposition of $\tau$. Let
$H=h(\bar{0})\oplus\sum _{i=1}^n \mathbb{C}K_i\oplus \sum _{i=0}^{n}\mathbb{C}d_i$
be our Cartan subalgebra for the root space decomposition of $\tau$. Define $\delta
_i,w_i \in H^{*}$ by setting 
\begin{center}
$\delta_i(h(\bar{0}))=0,\delta_i(K_j)=0$ and $\delta_i(d_j)=\delta_{ij}$. and
\end{center}
\begin{center}
$w_i(h(\bar{0}))=0$,$w_i(K_j)=\delta_{ij}$ and $w_i(d_j)=0$.
\end{center}
Take $\delta_{\beta}=\sum_{i=1}^n \beta_i\delta_i$ for $\beta \in \mathbb{C}^n$. For
$k\in \mathbb{Z}^n$, we shall refer to the vector $\delta_{k+\gamma}$ as the
translate of $\delta_k$ by the vector $\delta_{\gamma}$ where $\gamma \in
\mathbb{C}^n$. Define $g(\bar{k},\alpha):=\{x\in g(\bar{k})|[h,x]=\alpha
(h)x,\forall h\in h(\bar{0})\}$. Then  we have $\tau =\oplus_{\beta\in
\bigtriangleup }\tau_{\beta}$, where $\bigtriangleup \subseteq
\{\alpha+\delta_k|\alpha\in \bigtriangleup_{0,en},k\in \mathbb{Z}^n\}$. We have
$\tau_{\alpha+\delta_k}=g(\bar{k},\alpha)\otimes t^k$ for $\alpha\neq 0$, and
$\tau_{\delta_k}=g(\bar{k},0)\otimes t^k\oplus
\bigoplus_{i=1}^n\mathbb{C}t^kK_i\oplus \bigoplus_{i=1}^n\mathbb{C}t^kd_i $ for
$0\neq k \in \mathbb{Z}^n  $ and $\tau_0=H$.\vspace{5mm}
\\
In order to get a non-degenerate form on $H^*$ we extend $\alpha \in h(\bar{0})^*$
to $H$ by defining $\alpha(K_i)=\alpha(d_i)=0,\forall 1\leq i\leq n$. Then
$(h(\bar{0})|K_i)=0$, $(\delta_k|\delta_l)=(w_i|w_j)=0$ and
$(\delta_i|w_j)=\delta_{ij}$ and form on $h(\bar{0})$ is the restriction of the form
of $g$. One can easily check that this form is non-degenerate on $H^*$. A root
$\beta=\alpha+\delta_k$ is called a real root if $\alpha\neq 0$. Let $\bigtriangleup
^{re}$ denote the set of all real roots and
$\beta^{\vee}=\alpha^{\vee}+\frac{2}{(\alpha|\alpha)}\sum_{i=1}^n k_iK_i$ is a
co-root of $\beta$, where $\alpha^{\vee}$ is the co-root of $\alpha\in
\bigtriangleup_{0,en}$. For $\gamma\in \bigtriangleup^{re}$ define
$r_{\gamma}(\lambda)=\lambda-\lambda(\gamma^{\vee})\gamma$ for $\lambda\in H^*$. Let
$W$ be the Weyl group of $\tau$ generated by $r_{\gamma}, \forall \gamma\in
\bigtriangleup^{re}$. 
\begin{Def}
A $\tau$ -module $V$ is called integrable if 
\begin{itemize}
\item[(1)] $V=\bigoplus_{\lambda\in H^*}V_{\lambda}$, where $V_{\lambda}=\{v\in
V|h.v=\lambda(h)v \forall h\in H\}$ and $dim(V_{\lambda})\lneq \infty$
\item[(2)] All the real root vectors act locally nilpotently on $V,$ i.e.,
$g(\bar{k},\alpha)\otimes t^k$ acts locally nilpotently on $V$ for all $0\neq \alpha\in
\bigtriangleup_{0,en}$.
\end{itemize}
\end{Def}
Now we have the standard proposition 
\begin{Prop}\label{prop 3.2}
Let $V$ be an irreducible integrable module for $\tau$. Then
\begin{itemize}
\item[(1)] $P(V)=\{\gamma \in H^*|V_{\gamma}\neq 0\}$ is $W$- invariant.
\item[(2)] $dim(V_{\gamma})=dim (V_{w\gamma}), \forall w\in W$.
\item[(3)] If $\lambda \in P(V)$ and $\gamma\in \bigtriangleup^{re}$, then
$\lambda(\gamma^{\vee})\in \mathbb{Z}$.
\item[(4)] If $\lambda \in P(V)$ and $\gamma\in \bigtriangleup^{re}$, and
$\lambda(\gamma^{\vee})> 0$, then $\lambda-\gamma\in P(V)$.
\end{itemize}
\end{Prop}
 Now consider the natural triangular decomposition of $\tau$
 \begin{center}
 \[\tau^+=\bigoplus_{\substack{\alpha>0,k\in \mathbb{Z}^n}}g(\bar{k},\alpha)\otimes
t^k\]
 \end{center}
\begin{center}
 \[\tau^-=\bigoplus_{\substack{\alpha<0,k\in \mathbb{Z}^n}}g(\bar{k},\alpha)\otimes
t^k\]
 \end{center}
 \begin{center}
 \[\tau^0=\bigoplus_{k\in \mathbb{Z}^n}g(\bar{k},0)\otimes t^k\oplus Z\oplus S_n(m)\]
 \end{center}
Where $Z=\Omega_{A_n}/d_{A_n}$. 
\\
We have $V$ is an integrable module for $\tau$, so it will be integrable module for
$\widetilde{LT}=LT\bigoplus \Omega_{A_n}/d_{A_n}\bigoplus D$, where $D$ is the
$\mathbb{C}$ linear span of $d_i,1\leq i\leq n$.

\begin{lemma}\label{lem 3.1} 
Let $V$ be an integrable $\widetilde{LT}$-module . Then there exists some $\lambda
\in  P(V)$ such that $\lambda +\eta \notin P(V)$
\end{lemma}
 This is the lemma 3.6 of \cite{[4]}.
\begin{Th}
The highest weight space $V_{+}=\{v\in V |\tau^+ .v=0\}$ is non-zero.
\end{Th}
\begin{proof}
$\widetilde{LT}$ is a Lie subalgebra of $\tau.$ So an integrable module for $\tau$
with finite dimensional weight spaces  will be an integrable module for
$\widetilde{LT}$ with finite dimensional spaces. Now take $\lambda \in P(V)$ as in
lemma \ref{lem 3.1}. Using Proposition \ref{prop 3.2}, it can be proved that
$\lambda |_{h(\bar{0})}$ is a dominant integral weight of $g(\bar{0})$.Now  using
the same idea as \cite{[14]} (Theorem 2.4(ii)), we can find $\mu \in P(V)$ such that  $V_{\mu
+\alpha+\delta_k}=0$, $\forall \alpha\in \bigtriangleup_{0}^{+},k\in \mathbb{Z}^n$.
Again using Proposition \ref{prop 3.2} we will get  $V_{\mu +2\alpha+\delta_k}=0$,
$\forall \alpha\in \bigtriangleup_{0}^{+},k\in \mathbb{Z}^n$. Now our theorem is a
direct consequence of $\bigtriangleup=\bigtriangleup_{0,en}$.
\end{proof}
We can see that $V_+$ is a $\tau^0$ module.
\begin{lemma}\label{lem 3.2}

 (1) $V_+$ is irreducible module over $\tau^0$.\\
 \hspace*{23mm}
 (2) $V=U(\tau^-)V_+$.

\end{lemma}
\begin{proof}
Using PBW theorem and some weight arguments we can prove (1). Using (1) and PBW
theorem we can prove (2)
\end{proof}
\begin{lemma}\label{lem 3.3}
\begin{itemize}
\item[(1)] There exists uniqe  $\bar{\lambda}\in h(\bar{0})^*$ and a  $\beta \in
\mathbb{C}^n$, here $\beta$ need not be unique such that the weights of $V_+$ will
be of the form $\bar{\lambda}+\delta_{r+\beta}$ where $r\in \mathbb{Z}^n$ . 
\item[(2)] $\bar{\lambda}$ is a dominant integral.
\end{itemize}
\end{lemma}
\begin{proof}
(1) We can see that $h(\bar{0})$ commutes with $\tau^0$. Now using Lemma \ref{lem
3.2}  we will get $h(\bar{0})$ acts by scalars on $V_+$ and hence it will act by a
single linear functional on $V_+$. Let us denote it by $\bar{\lambda}$. If $D(u,r)
\in S_n(m)$ , then $D(u,r).V_{\mu}\subseteq V_{\mu+\delta_r}$ for $\mu \in
P(V)$.\\
(2) Direct consequence of Proposition \ref{prop 3.2} (4).
\end{proof}
If $\bar{\lambda}=0$ in Lemma \ref{lem 3.3}, then using Proposition \ref{prop 3.2},
we can check that only possible weights of $V$ are $\delta_{\beta},\beta\in
\mathbb{C}^n$. Now one can prove that $LT+\Omega_{A_n(m)}/d_{A_n(m)}$ acts trivially on
$V$. Then $V$ will be irreducible  module for $S_n(m)$ with finite dimensional
weight spaces with respect to the cartan $D=\oplus_{i=1}^n\mathbb{C}d_i$. It is
still unknown about the irreducible modules for $S_n(m)$ with finite-dimensional
weight spaces. So we will take $\bar{\lambda}\neq 0$ in Lemma \ref{lem 3.3}.\\
Take $\lambda \in P(V)$ from Lemma \ref{lem 3.1}. Then by our assumption there
exists $\alpha \in \bigtriangleup_0$ such that $\lambda({h_{\alpha}})\neq 0$. Then
using Lemma \ref{lem 3.2} we can get \[ V_+=\oplus _{r\in \mathbb{Z}^n}V_+(r),\]
where $V_+(r)=\{v\in V_+|d_i.v=(\lambda(d_i)+m_i)v,1\leq i\leq n\}$. So $V_+$ is
$\mathbb{Z}^n$ graded. Now for $r \in \Gamma$ we define
$V^{\prime}(r)=\oplus_{r_i\leq k_i< m_i+r_i}V_k$. Then $V_+=\oplus_{r \in
\Gamma}V_+^{\prime}(r)$ is $\Gamma$ graded.
\section{Action of Center}\label{sec 4}
In this section we will prove that $\Omega_{A_n(m)}/d_{A_n(m)}$ acts trivially on $V$. The following two lemmas are proved in \cite{[2]} as well as in \cite{[7]}

\begin{Prop}\label{prop 4.1}
There are only finitely many $h(\bar{0})\otimes A_n(m)\oplus
\Omega_{A_n(m)}/d_{A_n(m)}$- submodules of $V_+$
\end{Prop}
\begin{Prop}
The weight spaces of $V_+$ are uniformly bounded.
\end{Prop}
\begin{Prop}
$\Omega_{A_n(m)}/d_{A_n(m)}$ acts trivially on $V$.
\end{Prop}
\begin{proof}
From Proposition \ref{prop 4.1}, take a minimal $h(\bar{0})\otimes A_n(m)\oplus
\Omega_{A_n(m)}/d_{A_n(m)}$ -submodule of $V_+$, assume $W$. Now clearly $W$ will be
irreducible. Take $W^\prime=U(g(\bar{0})\otimes A_n(m)\oplus
\Omega_{A_n(m)}/d_{A_n(m)}\oplus D)W$. We can see that $W$ goes injectively to the
irreducible quotient of $W^{\prime}$. Now using \cite{[5]} we can see that
$\Omega_{A_n(m)}/d_{A_n(m)}$ acts trivially on $W$.One can easily check that  $S=\{v\in
V|t^rK_i.v=0, 1\leq i\leq n\}$ is a non-zero $\tau$-submodule of $V$. Hence we are
done by irreducibility of $V$.  
\end{proof}

\section{Action of $h_\alpha \otimes t^k$ on highest wieght space}\label{sec 5}
\begin{lemma}\label{lem 5.1}
Let $h\otimes t^k \in \tau^0$, where $k\in \Gamma \setminus \{0\}$ and $h \in
h(\bar{0})$. If there exists a nonzero vector of $V_+$ on which $h\otimes t^k $ acts
trivially, then  $h \otimes t^k $ is locally nilpotent on $V_+$.
\end{lemma}
\begin{proof}
Consider the set $S=\{v\in V_+|$ there is a $  n \in \mathbb{N}$ s.t. $(h\otimes
t^k)^n.v=0 \}$. By assumption clearly this set is nonzero.  We know that  $h\otimes
t^k$ commutes with $\bigoplus_{k\in \mathbb{Z}^n}g(\bar{k},0)\otimes t^k$. So
$(\bigoplus_{k\in \mathbb{Z}^n}g(\bar{k},0)\otimes t^k)S\subseteq S$. Now using the
equation 
\begin{center}
$(h\otimes t^k )^n.D(u,r)=D(u,r).(h\otimes t^k )^n-n(u.k)(h\otimes t^k
)^{n-1}.h\otimes t^{k+r}$
\end{center}
we can prove that $S_n(m).S\subseteq S$. Therefore $S$ is a nonzero $\tau^0$-
submodule of $V_+$. Hence by irreducibility of $V_+$ we can say that $S=V_+$. So the
result.  
\end{proof}
The idea of the proof of following two lemma is from \cite{[12]}, but the
calculations will be tottaly different. So we will give the proofs.
\begin{lemma}\label{lem 5.2}
Let $h\otimes t^k $ acts locally nilpotently on $V_+$, for every $k \in \Gamma
\setminus \{0\}$, where $h\in h(\bar{0})$. Then $h\otimes t^k $ acts trivially on
$V_+$, for every $k\in \Gamma\setminus \{0\}$.
\end{lemma}
\begin{proof}
Let $p \in \Gamma \setminus \{0\}$. Then by above discussion we will get $(h
\otimes t^{-p}.h \otimes t^p)^N=0$. Therefore $(h \otimes
t^{-p})^N.(h \otimes t^p)^N=0$.\\
Let $r_1\in \Gamma \setminus \{0\}$ be such that $r_1\notin \mathbb{Q}.p$. Then
there exists $u \in \mathbb{C}^n$ such that $(u,p)\neq 0$ and $(u,r_1)=0$. Therefore
we have the following
\begin{center}
$D(u,r_1).(h\otimes t^{-p})^N.(h \otimes t^p)^N=0$
\end{center}
This equation will give us \\
 
$N(u,p)h\otimes t^{r_1+p}.(h \otimes t^{-p})^N.(h\otimes
t^p)^{N-1}-N(u,p)h \otimes t^{r_1-p}.(h \otimes
t^{-p})^{N-1}.(h \otimes t^p)^{N-1}=0$\\
Now applying $h \otimes t^p$ on the above equation we will get
\begin{center}
$h \otimes t^{r_1-p}.(h \otimes t^{-p})^{N-1}.(h\otimes
t^p)^{N+1}=0$.
\end{center}
 Using induction on $j$, we want to prove, $h \otimes t^{r_1-p}.\dots
.h \otimes t^{r_j-p}.(h \otimes t^{-p})^{N-j}.(h \otimes
t^p)^{N+j}=0$.\\
 Where $1\leq j \leq N$ and $r_1,\dots, r_j \in \Gamma $ with $\sum_{i=1}^j
\epsilon_i r_i \notin \mathbb{Q}.p$ and $\epsilon_i \in \{0,1\}$.\\
 For $j=1$, we already prove this. Now assume the result holds for $1\leq j <N$.\\
 Let $r\in \Gamma$ with $r+ \sum_{i=1}^j \epsilon_ir_i \notin \mathbb{Q}.p$. So we
can find  $u \in \mathbb{C}^n$ such that $(u,p)\neq 0$ and $(u,r)=0$. Therefore we
have 
 \begin{center}
 $D(u,r).h \otimes t^{r_1-p}.\dots .h \otimes t^{r_j-p}.(h
\otimes t^{-p})^{N-j}.(h \otimes t^p)^{N+j}=0$.
 \end{center}
 From this equation we will get,
 \begin{center}
 $(N+j)(u,p)h \otimes t^{r_1-p}.\dots .h \otimes t^{r_j-p}.h
\otimes t^{r+p}.(h \otimes t^{-p})^{N-j}.(h \otimes
t^p)^{N+j-1}-(N-j)(u,p)h \otimes t^{r_1-p}.\dots .h\otimes
t^{r_j-p}.h \otimes t^{r-p}.(h \otimes t^{-p})^{N-j-1}.(h
\otimes t^p)^{N+j}=0$.
 \end{center}
 Applying $h\otimes t^p$, we have
 \begin{center}
 $h \otimes t^{r_1-p}.\dots .h \otimes t^{r_j-p}.h \otimes
t^{r-p}.(h\otimes t^{-p})^{N-j}.(h \otimes t^p)^{N+j}=0$.
 \end{center}
 Now taking $j=N$,
 \begin{center}
 $h \otimes t^{r_1-p}.\dots .h \otimes t^{r_N-p}.(h \otimes
t^p)^{2N}=0$
 \end{center}
 From this equation, using induction again, we can prove 
 \begin{center}
 $h \otimes t^{r_1-p}.\dots .h \otimes t^{r_{N+j}-p}.(h
\otimes t^p)^{2N-j}=0$,
  \end{center}
  for all $0\leq j\leq 2N$ with $\sum_{i=1}^{N+j}\epsilon_ir_i \notin \mathbb{Q}.p$
and $\epsilon_i \in \{0,1\}$.\\
  In this case instead of taking $D(u,r)$ we take $D(u,r-2p)$, it is easy to find
such $u \in \mathbb{C}^n$. Hence we will have
  \begin{center}
  $h \otimes t^{r_1-p}.\dots .h \otimes t^{r_{3N}-p}=0$,
  \end{center}
 with all $r_i \in \Gamma$ and $\sum_{i=1}^{3N}\epsilon_ir_i \notin \mathbb{Q}.p$. \\
 Let $s_1,\dots s_{3N}\in \Gamma \setminus \{0\}$ be arbitrary. Take $p \in \Gamma
\setminus \cup_{\epsilon_i \in \{0,1\}} \mathbb{Q}.(\sum_{i=1}^{3N}\epsilon_is_i)$
. Now taking $r_i=s_i+p$ in the above equation we can get 
 \begin{center}
 $h \otimes t^{s_1}\dots h \otimes t^{s_{3N}}=0.\cdots (*)$ .
 \end{center}
 Take $W=\{v \in V_+/ h_\alpha \otimes t^p.v=0, \forall p\in
\Gamma\setminus\{0\}\}$. By $(*)$, $W$ is nonzero. We can easily check that $W$ is
a $\tau^0$-submodule of $V_+$. Hence by the irreducibility of $V_+$, we have the
lemma.
 \end{proof}
\begin{lemma}\label{lem 5.3}
$h\otimes t^k$ acts injectively or trivially on $V_+$, for all $k\in \Gamma\setminus
\{0\}$ .
\end{lemma}
\begin{proof}
Let $h \otimes t^p$ acts injectively on $V_+$, for some $p \in \Gamma
\setminus \{0\}$.\\
\textbf{Claim:(1)} For any $q \in \Gamma \setminus \mathbb{Q}.p$, if $h
\otimes t^q$ acts locally nilpotently, then it will act nilpotently.\\
Assume $h \otimes t^q$ acts locally nilpotently. Since $dim(V_+(r))\leq N$
for all $r \in \Gamma$, we will have $(h \otimes t^{-q})^N.(h\otimes
t^q)^N=0$. \\
Suppose there exists $1\leq k \leq N$ and $N_k$, such that $(h \otimes
t^{-q})^k.(h \otimes t^q)^{N_k}=0$.\\ We assume $q \notin \mathbb{Q}.p$. So
there exists $u \i n \mathbb{C}^n$ such that $(u,p+q)=0$ and $(u,q)\neq 0$.\\
Now applying $D(u,p+q)$, to the above equation we will have\\ 
$N_k(u,q)h\otimes t^{p+q}.(h\otimes t^{-p})^k(h \otimes
t^{p})^{N_k-1}-k(u,q)h \otimes t^p. (h \otimes t^{-p})^{k-1}.(h
\otimes t^p)^{N_k}=0$.
Now applying $h\otimes t^q$, we have $h \otimes t^p. (h
\otimes t^{-p})^{k-1}.(h \otimes t^p)^{N_k+1}=0$. Since $h\otimes
t^p$ acts injectively on $V_+$, so $(h \otimes t^{-p})^{k-1}.(h
\otimes t^p)^{N_k+1}=0$. Now repeating the process we can find $N_0\in \mathbb{N}$
such that $(h_\alpha \otimes t^q)^{N_0}=0$.\\
\textbf{Claim:(2)} If $q\in \Gamma \setminus \{0\} \cap \mathbb{Q}.p $, then
$h \otimes t^q$ acts injectively on $V_+$.\\
First we will prove that $h \otimes t^{\frac{a}{b}p}$ acts injetively, if
$a\in \mathbb{Z}_-$, $b \in Z_+$.\\
If for some such $a,b \in \mathbb{Z}$, $h \otimes t^{\frac{a}{b}p}$ is not
acting injectively  then by lemma \ref{lem 5.1}, it will act locally nilpotently.
Therefore $(h \otimes t^p)^{-a}.(h \otimes t^{\frac{a}{b}p})^b$ acts
nilpotently on each homogeneous spaces. Hence $(h \otimes
t^p)^{-aN}.(h \otimes t^{\frac{a}{b}p})^{bN}=0$. Therefore $(h \otimes
t^{\frac{a}{b}p})^{bN}=0$, since $h \otimes t^p$ acts injectively on $V_+$.
Assume $N_0 \in \mathbb{N}$ be smallest such that $(h \otimes
t^{\frac{a}{b}p})^{N_0}=0$. Let $q \in \Gamma \setminus \mathbb{Q}.p$, then we can
find $u \in \mathbb{C}^n$ such that, $(u,q-\frac{a}{b}p)=0$ and $(u,p)\neq 0$. Now
take 
\begin{center}
$0=D(u,q-\frac{a}{b}p).(h \otimes
t^{\frac{a}{b}p})^{N_0}=N_0(u,\frac{a}{b}p)(h\otimes
t^{\frac{a}{b}p})^{N_0-1}. h \otimes t^q. V_+$
\end{center}
 From the minimality of $N_0$, we can say from the above equation and lemma \ref{lem
5.1} that $h \otimes t^q$ acts locally nilpotently. Therefore by claim (1)
$h \otimes t^q$ acts nilpotently on $V_+$. Choose $N^\prime$ be the smallest
s.t, $(h \otimes t^q)^{N^\prime}=0$. Choose $u \in \mathbb{C}^n$ s.t
$(u,p-q)=0$ and $(u,q)\neq 0$. Then $0=D(u,p-q)(h\otimes
t^q)^{N^\prime}=N^\prime (u,q)(h \otimes t^q)^{N^\prime-1}.h \otimes
t^p$.\\ Now by minimality of $N^\prime$, we can see that $h\otimes t^ p$ is
not acting injectively, which is a contradiction. So our claim (2) is true for $a
\in \mathbb{Z}_-,b \in \mathbb{Z_+}$. Now let both $a,b$ are positive. We know by
above $h \otimes t^{-ap}$ acts injectively. By our assumption $h
\otimes t^{-ap}. (h \otimes t^{\frac{a}{b}p})^b$ acts nilpotently on each
homogeneous spaces. Now repeating the same process as above we will get the
contradiction, $h \otimes t^p$ not acting injectively on $V_+$. So our claim
(2) is true. We are now going to prove that $h\otimes t^q$ acts injectively
for all $q \in \Gamma \setminus \{0\}$. We already prove it for $q \in
\mathbb{Q}.p$. Let $q \notin \mathbb{Q}.p$ and $h \otimes t^q$ does not acts
injectively on $V_+$.  Therefore by lemma \ref{lem 5.1} and claim(1), we can say
that $h\otimes t^q$ acts nilpotently on $V_+$. Using the minimal power of
$h \otimes t^q$, for which it acts trivially,  we will arrive at $h
\otimes t^p$ not acting injectively on $V_+$, a contradiction to our assumption.
Hence using lemma \ref{lem 5.1} and \ref{lem 5.2} we have the lemma. 

\end{proof}

Now we have choosen $\alpha\in \bigtriangleup_0$ such that
$\bar{\lambda}(h_{\alpha})\neq 0$. Now if there is a $k\in \Gamma\setminus \{0\}$
such that $h_{\alpha}\otimes t^k$ is not injective on $V_+$, then by Lemma \ref{lem
5.3} we can say that $h_{\alpha}\otimes t^m$ acts trivially on $V_+$. Then using
\cite{[6]} we can say that $h_{\alpha}$ acts trivially on $V_+$, which is a
contradiction to our assumption.\\
We know $V_+=\oplus_{r\in \Gamma}V_+(r)$. Let $\{v_0,\dots ,v_k\}$ be a basis of
$V_+(0)$. Therefore by above discussion, $\{h_\alpha \otimes t^p.v_0, \dots,
h_\alpha \otimes t^p.v_k \}$ will be a basis of $V_+(p)$ for all $p \in \Gamma$. Now
assume, $h_\alpha \otimes t^p.h_\alpha \otimes t^q.(v_1, \dots , v_k)=h_\alpha
\otimes t^{p+q}.(v_1,\dots , v_k)A_{p,q}$, where $A_{p,q}\in M_k(\mathbb{C})$. Since
$h_\alpha \otimes t^p$ are commuting family, we can easily see that $\{A_{p,q}| p,q
\in \Gamma\}$ is a commuting family. Therefore we can find an invertible matrix $S
\in M_k(\mathbb{C})$ such that $\{B_{p,q}=S^{-1}A_{p,q}S|p,q \in \Gamma\}$ is upper
triangular family. Set $(w_1, \dots ,w_k)=(v_1,\dots , v_k).S$, then we have
$h_\alpha \otimes t^p.h_\alpha \otimes t^q.(w_1, \dots , w_k)=h_\alpha \otimes
t^{p+q}.(w_1,\dots , w_k)B_{p,q}$. So without loss of generality we can assume that
$A_{p,q}$ is upper triangular for all $p, q \in \Gamma$.
\begin{lemma}\label{lem 5.4}
Let $p,q \in \Gamma$, then there exists $\lambda_{p,q}$ and $v_0 \in V_+(0) $ such
that $(h_\alpha \otimes t^p.h_\alpha \otimes t^q-\lambda_{p,q}h_\alpha \otimes
t^{p+q})v_0=T_{p,q}v_0=0$.
\end{lemma}
\begin{proof}
Since for $\lambda \in \mathbb{C}$, $(h_\alpha \otimes t^p.h_\alpha \otimes
t^q-\lambda h_\alpha \otimes t^{p+q}).(v_1,\dots v_k)=h_\alpha \otimes
t^{p+q}.(v_1,\dots, v_k)(A_{p,q}-\lambda . I_k)$. Now the Lemma follows from
$A_{p,q}$ has eigenvalues.
\end{proof}
\begin{lemma}\label{lem 5.5}
For $p,q \in \Gamma$, $\lambda \in \mathbb{C}$, if there exists $v \in
V_+(0)\setminus \{0\}$ such that $T_\lambda .v=0$, then $T_\lambda$ acts locally
nilpotently on $V_+$.
\end{lemma}
\begin{proof}
The proof will be same as Lemma \ref{lem 5.1}.
\end{proof}
Now using Lemma \ref{lem 5.4} and Lemma \ref{lem 5.5}, we can say that for every $p,q \in
\Gamma$, $A_{p,q}$ has exactly one eigenvalue (say, $\lambda_{p,q}$).
\begin{Prop}\label{prop 5.1}
\itemize 
\item[(1)] $h_{\alpha}\otimes t^k$ acts injectively on $V_+$ for every $k \in \Gamma$.
\item[(2)] $h_{\alpha}\otimes t^r.h_{\alpha}\otimes
t^s=\lambda_{r,s}h_{\alpha}\otimes t^{r+s}$ on $V_+$, where $\lambda_{r,s}=\lambda$
for all $r\neq 0,s\neq 0,r+s\neq 0$, $\lambda_{r,-r}=\mu$ for all $r\neq 0$ and
$\lambda_{0,r}=\bar{\lambda}(h_{\alpha})$ for all $r\in \Gamma$. Further we have
$\mu \lambda_{0,r}=\lambda^2\neq 0$.
\item[(3)] dim $ (V_+^{\prime}(r))=$ dim $ (V_+^{\prime}(s))$ for all $r,s \in
\Gamma$. 
\end{Prop} 
\begin{proof}
(1) is a direct consequence of above discussion. (3) also follows from (1). Now we
will prove (2). Set $V_0=\sum _{i=0}^{k-1}\mathbb{C} v_i$, where $v_0=0$. Since
$A_{p,q}$ is upper triangular, we have $T_{p,q}.v_l \in \sum_{i=0}^{l-1}h_\alpha
\otimes t^{p+q}.v_i$. We can easily see that $T_{p,q}.V_+(-p-q) \subseteq V_0$. Now
for any $p,q,r \in \Gamma$,with $u \in \mathbb{C}^n$ such that $(u,r)=0$, define
$S_{p,q}^r=(u,p)h_\alpha \otimes t^{p+r} .h_\alpha \otimes t^q+(u,q)h_\alpha \otimes
t^p. h_\alpha \otimes t^{q+r}-\lambda_{p,q}(u,p+q)h_\alpha \otimes t^{p+q+r}$.
Therefore $h_\alpha \otimes
t^{-(p+q+r)}.S_{p,q}^r.v_k=\lambda_{-(p+q+r),p+q+r}((u,p)\lambda_{p+r,q}+(u,q)\lambda_{p,q+r}-(u,p+q)\lambda_{p,q}).v_k+w^\prime$,
where $w^\prime \in V_0$. One can also check that $h_\alpha \otimes
t^{-(p+q+r)}.S_{p,q}^r.V_0 \subseteq V_0$. Now we can see that
$[D(u,r),T_{p,q}]=S_{p,q}^r$ . Given any $p,q,r \in \Gamma$, and $u \in
\mathbb{C}^n$ with $(u,r)=0$ we have\\

$0=(h _\alpha \otimes t^{-(p+q+r)})^k.(D(u,r))^k.(T_{p,q})^k.v_k=T_{p,q}.w+k!(h
_\alpha \otimes t^{-(p+q+r)})^k.(S_{p,q}^r)^k.v_k$, where $w \in V_+(-p-q)$.\\
From the above equation, we will have
\begin{center}
$(u,p)\lambda_{p+r,q}+(u,q)\lambda_{p,q+r}-(u,p+q)\lambda_{p,q}=0, \cdots (**)$,
\end{center}
for any $p,q,r \in \Gamma$.\\
Take $W=\{v \in V_+|T_{p,q}.v=0\}$, then $v_1 \in W$. Hence using the above equation
we can check that $W$ is a nonzero $\tau^0$-submodule of $V_+$, therefore by
irreducibility of $V_+$, $W=V_+$.\\
Now taking $p=jq$, with $j\mathbb{Q}\setminus \{0,-1\}$ in $(**)$, we will get
\begin{center}
$j(u,q)\lambda_{jq+r,q}+(u,q)\lambda_{jq,q+r}-(j+1)(u,q)\lambda_{p,q}=0$.
\end{center}
Taking, $r \in \Gamma \setminus \mathbb{Q}.q$, we can see that 
\begin{center}
$j\lambda_{jq+r,q}+\lambda_{jq,q+r}-(j+1)\lambda_{q,q}=0$.
\end{center}
Putting $j=1$, we have $\lambda_{q,q+r}=\lambda_{q,q}$ and hence
$\lambda_{jq+r,q}=\lambda_{jq,q+r}=\lambda_{q+r,q+r}=\lambda_{q,q}$ for all $r \in
\Gamma \setminus \mathbb{Q}.q$. Therefore $\lambda_{p,q} =\lambda_{q,q}(=\lambda$, say), if $p\neq
0,q\neq 0,p+q\neq 0$.\\
Taking $r=-p-q$ with $p \in \Gamma \setminus \mathbb{Q}.q$ in $(**)$, we have
$\lambda_{p,-p}=\lambda_{q,-q}(=\mu$,say) and hence for every $p,q \in \Gamma$. It is easy to
see that $\lambda_{r,0}=\bar{\lambda}(h_\alpha)$.\\
Let us assume $p,q \in \Gamma$ be such that $p\neq 0,q\neq 0,p+q\neq 0$. then for any $v \in V_+$, we have \\
\begin{center}
$(h_\alpha \otimes t^p.h_\alpha \otimes t^{-p}).(h_\alpha \otimes t^q. h_\alpha \otimes t^{-q}).v=(\mu)^2.(\bar{\lambda}(h_\alpha))^2v=\lambda_{p,q} \lambda_{-p,-q}\lambda_{p+q,-p-q}\bar{\lambda}(h_\alpha).v$.
\end{center}
hence we have the identity of (2).
\end{proof}
\section{classification}\label{sec 6}
Now recall that our Lie algebra reduces to $\tau= LT\oplus S_n(m)$ with
$\tau^0=\oplus _{k\in \mathbb{Z}^n}g(\bar{k},0)\otimes t^k\oplus S_n(m)$. Take
$g^{\prime}=\{x\in g |[h,x]=0,\: \forall\:h \in h(\bar{0})\}$. Now since
$g^{\prime}$ is invariant under $\sigma_i$'s, therefore $g^{\prime}$ is $\Lambda$
graded. It is easy to see that $L(g^{\prime},\sigma)=LT^0$. Let us take
$S_n^{\prime}(m)=$ span$\{D(u,r)-D(u,0)|u\in \mathbb{C}^n,r \in \Gamma,(u.r)=0\}$.
We can easily check that $S_n^\prime (m)$ is a Lie subalgebra of $S_n(m)$.
Furthermore let us set $L=S_n^\prime(m) \ltimes L(g^\prime,\sigma)$ and $W=$
span$\{h_\alpha \otimes t^r.v-v|r\in \Gamma,v \in V_+\}$. We can see that $W$ is a
$L$ module.
\begin{lemma}
\itemize
\item[(1)] $W$ is a proper $L$-submodule of $V_+$.
\item[(2)] $\widetilde{V}=V_+/W$ is a finite dimensional $L$-module.

\end{lemma}
\begin{proof}
The proof is parallel to Lemma 7.2 of \cite{[7]}. Let $z_i=h_{\alpha}\otimes t^m_i$
for each $i=1,\dots,n$. Without loss of generality we can assume that
$\lambda_{r,s}=1$ for $r\neq 0,s\neq 0,r+s\neq 0$. Therefore we can say that $W=$
span$\{z_i.v-v|v\in V_+\}$. Now the proof will be same as \cite{[7]}
\end{proof}
Let $\beta\in\mathbb{C}^n$ be as in Lemma \ref{lem 3.3}. Then for any $L$ module
$V^\prime$ we can give a $\tau^0$ module structure on $L(V^\prime)$ by 
\begin{center}
$x\otimes t^k.(v_1\otimes t^s)=((x\otimes t^k).v_1)\otimes t^{k+s}$.
\end{center}
\begin{center}
$D(u,r).(v_1\otimes t^s)=((D(u,r)-D(u,0)).v_1)\otimes t^{r+s}+
(u,s+\beta)(v_1\otimes t^{r+s})$
\end{center}
for all $v_1 \in V^\prime ,x \in g(\bar{k},0)$ and $D(u,r)\in S_n(m)$. \\
For $v\in V_+$, let $\bar{v}$ be the image of $v$ in $\widetilde{V}$. Now define 
\begin{center}
$\phi: V_+ \rightarrow L(\widetilde{V})$
\end{center}
\begin{center}
by $v\mapsto \bar{v}\otimes t^k$ for $v\in V_+(k).$
\end{center}
This map is clearly a nonzero $\tau^0$ -module homomorphism. Hence by irreducibility
of $V_+$ it follows that $V_+\cong \phi (V_+)$ is a $\tau ^0$ submodule of
$L(\widetilde {V})$. \\Clearly $L$ is naturally $\Lambda$ graded. Now since $V_+$
and $W$ is $Z^n$ graded, therefore they are naturally $\Lambda$ graded and hence so
$\widetilde{V}$. Therefore $\widetilde{V}=\oplus_{\bar{p}\in
\Lambda}\widetilde{V}(\bar{p})$.\\
Now for $\bar{p} \in \Lambda$ we set 
\begin{center}
$L(\widetilde{V})(\bar{p})=\{v\otimes t^{k+r+p}|v\in \widetilde{V}(\bar{k}),r \in
\Gamma,k\in \mathbb{Z}^n\}$
\end{center}
It can be easily verified that $L(\widetilde{V})(\bar{p})$ is a $\tau^0$ submodule
of $L(\widetilde{V})$.
The following result can be deduced similarly as in \cite{[3]} (Same result is also proved in \cite{[3]}).
\begin{Prop}\label{prop 6.1}
\itemize
\item[(1)] $V_+\cong L(\widetilde{V})(\bar{0})$ as $\tau^0$ -modules.
\item[(2)] $\widetilde{V}$ is $\Lambda$ -graded-irreducible module over $L$.
\item[(3)] $\widetilde{V}$ is completely reducible module over $L$ and all its
irreducible components are mutually isomorphic as $S_n^{\prime}(m)\ltimes
h(\bar{0})\otimes A(m)$ -modules.
\end{Prop}
Now we will concentrate on irreducible representation of $L$. Let $(W,\pi)$ be a
finite dimensional representation of $L$. Let $\pi (L(g^{\prime},\sigma))=g^1$, then
$\pi(L)=g^1\oplus g^2$, where $g^2$ is the unique complement of $g^1$ in $gl(W)$
(Proposition 19.1(b) of \cite{[8]} ). So $W$ will be an irreducible module for
$g^1\oplus g^2$. Therefore $W\cong W_1\otimes W_2$, where $W_1$ and $W_2$ are
irreducible modules for $g^1$ and $g^2$ respectively ( \cite{[9]} ). Let
$g^{\prime}=g^{\prime}_{ss}\oplus R$, where $g^{\prime}$ and $R$ are Levi and
radical part of $g^\prime$. Then as $\sigma_i(g^\prime)=g^\prime$ and
$\sigma_i(R)=R$ for $1\leq i\leq n$, we have
$L(g^\prime,\sigma)=L(g^\prime_{ss},\sigma)\oplus L(R,\sigma)$. Now $W_1$ is
irreducible module for $L(g^\prime,\sigma)$. As $R$ is solvable ideal, it follows
that $\pi(L(R,\sigma))$ lies in the center of $\pi(L)$, which is atmost one
dimensional. Hence $L(R,\sigma)$ acts as a scalar on $W$. So $W_1$ will be a
irreducible module for $L(g^\prime_{ss},\sigma)$. \\

Fix a positive integer $l$. For each $i$, let $a_i=(a_{i,1},\dots ,a_i,l)$ such that
$a_{i,j}^{m_i}\neq a_{i,t}^{m_i}$ for $j\neq t$ (*). Let $g$ be a finite dimensional
semisimple Lie algebra. Let $\sigma_1, \dots \sigma_n$ be finite order automorphisms
on $g$ of order $m_1,\dots m_n$ respectively. Let $L(g,\sigma)$ be the corresponding
multiloop algebra. Let $I=\{(i_1,i_2,\dots, i_n)|1\leq i_j\leq l\}$. Now for
$S=(i_1,i_2,\dots ,i_n)\in I$ and $r=(r_1,r_2,\dots r_n)\in \mathbb{Z}^n$,
$a_S^r=a_{1,i_1}^{r_1}a_{2,i_2}^{r_2}\cdots a_{n,i_n}^{r_n}.$ Now consider the
evaluation map $\phi : g\otimes A\rightarrow \bigoplus g$ ($l^n$ copies), $\phi
(X\otimes t^r)=(a_{I_1}^r,a_{I_2}^r, \dots,a_{I_{l^n}}^r)$, where $I_1,I_2, \cdots
I_{l^n}$ is some order on $I$. Now consider the restriction of $\phi$ to
$L(g,\sigma)$.
\begin{Th}
Let $W^\prime$ be a finite dimensional irreducible representation of $L(g,\sigma)$.
Then the representation factors through $\bigoplus g$ ($l^n$ copies).
\end{Th}\label{Th 6.2}
The proof of this result comes from \cite{[10]} .
\begin{Prop}
Let $W_1$ be irreducible module for $L(g^\prime _{ss},\sigma)$ as above. Then the
representation of $L(g^\prime _{ss},\sigma)$ factors through only one copy of
$\bigoplus g^\prime _{ss}$. So $g^1_{ss}\cong g^\prime_{ss}$.
\end{Prop}
\begin{proof}
We know by Theorem \ref{Th 6.2} that the representation factors through $l^n$
copies, for some positive integer $l$. we will prove here that $l=1$. The proof is
same as \cite{[2]}. Choose the $i$-th piece of $g^\prime_{ss}$ and choose the
proection of the map $\pi$, say $\pi_i$ onto it. Doing the same calculation as in
\cite{[2]} we will get $\pi_i(S_n^{\prime}(m))=0$ and $a_{I_i}^r=1$ for all $r\in
\Gamma$. Now suppose there is atleast two pieces, say $i$-th and $j$-th piece is
there. Therefore $I_i$ and $I_j$ are two different element of $I$ with
$a_{I_i}^r=1=a_{I_j}^r$ for all $r \in \Gamma$. Let $I_i=(i_1,i_2,\dots, i_n)$ and
$I_j=(j_1,j_2,\dots j_n)$. Therefore there is $k$ with $1\leq k\leq n$ such that
$i_k\neq j_k$. Now if we take $r=(0,\dots, m_k,\dots, 0)$ then
$a_{I_i}^r=1=a_{I_j}^r$ will give $a_{k,i_k}^{m_k}=a_{k,j_k}^{m_k}$, a contradiction
to (*). So there is atmost one piece.
\end{proof}
Now we know $\pi_i(S_n^\prime(m))=0$, therefore $g^2\subseteq \pi(S_n^\prime(m))$.
Now our aim is to understand  finite dimensional irreducible modules for
$S_n^\prime(m)$. There is a relation between finite dimensional $S_n^\prime$ modules
and $S_n\ltimes A$-modules with finite dimensional weight spaces.

\begin{Th}
Suppose $W$ is a finite dimensional irreducible $sl_n$-module (extend $W$ to $gl_n$
by letting $I$ ats trivially). Let $\alpha,\beta \in \mathbb{C}^n$. Take
$L(W)=W\otimes A$ and the action $D(u,r)(w\otimes t^k)=(u.k+\beta)w\otimes
t^{k+r}+\sum_{i,j}u_ir_jE_{ji}w\otimes t^{k+r}$, for $r\neq 0$ and $D(u,0)(w\otimes
t^k)=(u.\alpha+k)w\otimes t^k$ and $t^r(w\otimes t^k)=w\otimes t^{k+r}$. And all
irreducible representations of $S_n\ltimes A$ with finite dimensional weight spaces
occur in this way.
\end{Th}
This result is proved in \cite{[11]}\\
\begin{Th}
Let $W$ be finite dimensional irreducible module for $sl_n$ and extend it trivially
as before. Let $E_{ij}$ be generators of $gl_n$. Then $W$ can be made into
$S_n^\prime$-module by the action: $(D(u,r)-D(u,0)).w=\sum_{i,j}E_{ji}w+(u.\zeta)w$,
where $\zeta \in \mathbb{C}^n$. In this case $W$ will be irreducible module for
$S_n^\prime$. Every irreducible finite dimensional $S_n^\prime$ module occur in this
way.

\end{Th}
\begin{proof}
See Theorem 4.5 of \cite{[2]} and discussion thereafter.
\end{proof}
We know $\widetilde{V}$ is completely reducible $L$ module. Therefore
$\widetilde{V}=\oplus_{i=1}^K \widetilde{V}_i$ for some $K\in \mathbb{N}$. Then by
the previous discussion each $\widetilde{V}_i\cong W_1^{i}\otimes W_2^i$ as
$g^\prime_{ss}\oplus sl_n$ module, where $W_1^i$, $W_2^i$ are irreducible modules
for $g^\prime_{ss}$ and $sl_n$ respectively. since each component $\widetilde{V}_i$
are isomorphic as $S_n^\prime  (m)\ltimes (h(\bar{0})\otimes A(m))$ modules, we can
take $W_2^i\cong W_2^1$ ($ =W_2$ ,say) as $sl_n$-modules for each $i \in \{1,\dots
,K\}$. Now consider $W_1=\sum_{i=1}^KW_1^i$, which is a
$L(g^\prime_{ss},\sigma)$-module,  in particular $g^\prime_{ss}$ module. Since each
$W_1^i$ are irreducible, without loss of generality we can assume that the above sum
is direct. It is easy to see that $L$ is $\Lambda$-graded with zero-th component
$S_n^\prime  (m)\ltimes (h(\bar{0})\otimes A(m))$ and since $\widetilde{V}$ is
$\Lambda$ graded irreducible module (\ref{prop 6.1}), we can take $W_1$ as
$\Lambda$-graded irreducible $L(g^\prime_{ss},\sigma)$ module and $W_2$ is zero
graded as $S_n^\prime  (m)$ lies inside the zero-th graded component of $L$.\\
We now define a $\tau^0$-module structure on $W_1\otimes W_2\otimes A_n$ by 
\begin{center}
$X\otimes t^k(w_1\otimes w_2\otimes t^l)=Xw_1\otimes w_2\otimes t^{k+l}$, for $k,l
\in \mathbb{Z}^n$ and $X\in g^\prime _{ss}(\bar{k})$.
\end{center}
\begin{center}
$D(u,r)(w_1\otimes w_2\otimes t^k)=(u.k+\beta)w_1\otimes w_2\otimes
t^{k+r}+w_1\otimes((\sum_{i,j})u_ir_jE_{ji}w_2)\otimes t^{k+r}$ for $r\neq 0$ and
$D(u,r)\in S_n(m)$.
\end{center}
\begin{center}
$D(u,0)(w_1\otimes w_2 \otimes t^k)=(u.k+\alpha)w_1\otimes w_2\otimes t^k$.
\end{center}
Now take any one dimensional representation of $L(R,\sigma)$ say $\psi$. Then for
$y\in R(\bar{k})$ we take $y\otimes t^k(w_1\otimes w_2\otimes
t^l)=\psi(y)(w_1\otimes w_2 \otimes t^{k+l})$.  Since $W_1$ is $\Lambda$ graded
which is compatible with $\Lambda$-gradation of $g^\prime_{ss}$, so the submodule
$V^\prime=\bigoplus_{k \in \mathbb{Z}^n} W_{1,k}\otimes W_2\otimes t^k$ will be
irreducible module for $\tau^0$. One can easily check that
$L(\widetilde{V})(\bar{0})\cong V^\prime$ as $\tau^0$-module.\\
Now as in \cite{[7]} consider the triangular dicomposition of $LT$ given by 
\begin{center}
$LT=\tau^{-}\oplus L(g^\prime_{ss},\sigma)\oplus L(R,\sigma)\oplus \tau^+$.
\end{center}
 Let $V_1$ be the unique $\Lambda$-graded irreducible quotient of the induced module
for $W_1$ and let $(V_1^i,\rho _i)$ denote the irreducible quotient of the induced
representation of $W_1^i$, for each $1\leq i \leq K$. Now by \cite{[7]} we can say
that $V_1=\oplus_{i=1}^KV_1^i$ is $\Lambda$-graded irreducible module for $LT$.\\
 Define $\tau$-module structure on $V_1\otimes W_2\otimes A$ by 
 \begin{center}
 $(x\otimes t^l).(\sum_{i=1}^K v_1^i\otimes v_2\otimes t^k)=\sum_{i=1}^K
(\rho_i(x\otimes t^l)v_1^i)\otimes v_2\otimes t^{k+l}$,
\end{center}  
\begin{center}
$D(u,r).(\sum_{i=1}^Kv_1^i\otimes v_2\otimes
t^k)=(u.k+\beta)(\sum_{i=1}^Kv_1^i\otimes v_2\otimes
t^{k+r})+\sum_{i=1}^k\sum_{l,j}v_1\otimes (u_lr_jm_jE_{jl}).v_2)\otimes t^{k+r}$.
\end{center}
Let $V_1=\oplus _{\bar{k}\in \Lambda}V_{1,\bar{k}}$. Take $\bigoplus_{k\in
\mathbb{Z}^n}V_{1,\bar{k}}\otimes W_2\otimes t^k$, which is a $\tau$-submodule of
$V_1\otimes W_2\otimes A$. One can prove that $\bigoplus_{k\in
\mathbb{Z}^n}V_{1,\bar{k}}\otimes W_2\otimes t^k$ is irreducible $\tau$ module.
\begin{Th}\label{Thm 6.6}
Let $V$ be an irreducible integrable $\tau$ module with finite dimensional weight
spaces, with all $K_i$ acting trivially for $1\leq i\leq n$. Then $V\cong
\bigoplus_{k\in \mathbb{Z}^n}V_{1,\bar{k}}\otimes W_2\otimes t^k$.
\end{Th}

 \end{document}